\documentclass[10pt,conference]{IEEEtran}
\newtheorem{theorem}{Theorem}
\newtheorem{lemma}{Lemma}

\newtheorem{corollary}{Corollary}
\newtheorem{conjecture}{Conjecture}
\begin{document}


\title{Reconstruction of permutations distorted by single transposition errors}

\author{
\authorblockN{Elena Konstantinova}
\authorblockA{Sobolev Institute of Mathematics \\
Siberian Branch of Russian Academy of Sciences \\
Novosibirsk, Russia \\
e\_konsta@math.nsc.ru} \and
\authorblockN{Vladimir Levenshtein}
\authorblockA{Keldysh Institute of Applied Mathematics \\
Russian Academy of Sciences \\
Moscow, Russia \\
leven@keldysh.ru}
\and
\authorblockN{Johannes Siemons}
\authorblockA{School of Mathematics\\University of East Anglia \\
Norwich, UK \\
j.siemons@uea.ac.uk}}

\maketitle

\begin{abstract}
The reconstruction problem for permutations on $n$ elements from
their erroneous patterns which are distorted by transpositions is
presented in this paper. It is shown that for any $n \geq 3$ an
unknown permutation is uniquely reconstructible from 4 distinct
permutations at transposition distance at most one from the unknown
permutation. The {\it transposition distance} between two
permutations is defined as the least number of transpositions needed
to transform one into the other. The proposed approach is based on
the investigation of structural properties of a corresponding Cayley
graph.  In the case of at most two transposition errors it is shown
that $\frac32(n-2)(n+1)$ erroneous patterns are required in order to
reconstruct an unknown permutation. Similar results are obtained for
two particular cases when permutations are distorted by given
transpositions. These results confirm some bounds for regular graphs
which are also presented in this paper.
\end{abstract}

\section{Introduction}
\noindent Efficient reconstruction of arbitrary sequences was introduced and
investigated by Levenshtein for combinatorial channels with errors
of interest in coding theory such as substitutions, transpositions,
deletions and insertions of symbols \cite{lev97, eros}. Sequences
are considered as elements of a vertex set $V$ of a graph $\Gamma=(V,E)$
where an edge $\{x,y\}\in E$ is viewed as the single error transforming
$x$ into $y\in V$.  One of the metric problems which arises here is
the problem of reconstructing an unknown vertex $x\in V$ from a
minimum number of vertices in the metric ball $B_r(x)$ of radius $r$
centered at the vertex $x\in V$. It is reduced to finding the value

\begin{equation}
\label{e1} N(\Gamma,r)=\max_{x,y\in V(\Gamma),\ x\neq y}| B_r(x)\cap
B_r(y)|,
\end{equation}
since $N(\Gamma,r)+1$ is the least number of distinct vertices in
the ball $B_r(x)$ around the unknown vertex $x$ which are sufficient
to reconstruct $x$ subject to the condition that at most $r$ single
errors have happened. As one can see, this problem is based on
considering metric balls in a graph but it differs from traditional
packing and covering problems in various ways. It is motivated by a
transmission model where  information is realized in the presence of
noise without encoding or redundancy,  and where the ability to
reconstruct a message (vertex) uniquely depends on having a
sufficiently large number of erroneous patterns of this message.

The value (\ref{e1}) was studied for the Hamming and Johnson graphs
\cite{eros}. Both graphs are distance--regular and the first is a Cayley graphs. The
problem of finding the value (\ref{e1}) is much more complicated for
graphs which are not distance--regular. Cayley graphs of this kind arise for instance on the symmetric group and the signed
permutation group, when the reconstruction of permutations and signed
permutations is considered for distortions  by single reversal
errors \cite{Kon05_1,Kon05_2}.

In this paper we continue these investigations and consider the
reconstruction problem for permutations distorted by single
transposition errors which consist of swapping\, 1) any two elements
of a permutation; 2) any two neighboring elements of a permutation;
and 3) the first and any other element of a permutation. The
corresponding graphs are the  transposition Cayley graph, the
bubble--sort Cayley graph and the star Cayley graph. They are
regular but not distance--regular. We investigate the combinatorial
properties of these graphs and present the values (\ref{e1}) when
$r=1,2$ in each case. Some bounds on $N(\Gamma, 1)$ and $N(\Gamma,
2)$ for regular graphs are also considered. It is shown that the
bubble--sort and star Cayley graphs are examples for which these
bounds are attained.

\section{Definitions, notation, general results}

\noindent Let $G$ be a finite group and let $S$ be a set of generators of $G$
such that the identity element $e$ of $G$ does not belong to $S$ and
such that $S=S^{-1},$ where $S^{-1}=\{s^{-1}: s \in S\}$. In the
{\it Cayley graph} $\Gamma=Cay(G,S)=(V,E)$ vertices correspond to
the elements of the group, i.e. $V=G$, and
 edges correspond to multiplication on the right by  generators, i.e.
$E=\{\{g,gs\}: g\in G, s\in S\}.$ Denote by $d(x,y)$ the {\it path
distance} between the vertices $x$ and $y$ in $\Gamma,$ and by
$d(\Gamma)={\rm max}\{d(x,y)\,:\,x,y\in V\}$ the {\it diameter} of
$\Gamma$. In other words, in a Cayley graph the diameter is the
maximum, over $g\in G,$ of the length of a shortest expression for
$g$ as a product of generators. For the vertex $x$ let
$S_r(x)=\{y\in V\, :\, d(x,y)=r \}$ and $B_r(x)=\{y\in V\,:\, d(x,y)
\leq r\}$ be the {\it sphere} and the {\it ball} of radius $r$
centered at $x,$ respectively. The vertices $y\in B_r(x)$ are {\it
$r$-neighbors} of the vertex $x$.

As mentioned in the Introduction, the value (\ref{e1}) was
investigated initially for distance--regular graphs such as the
Hamming and Johnson graphs. Let us recall that a simple connected
graph $\Gamma$ is  {\it distance--regular} if there are integers
$b_i, c_i$ for  $i\geq 0$ such that for any two vertices $x$ and $y$
at distance $d(x,y)=i$ there are precisely $c_i$ neighbors of $y$ in
$S_{i-1}(x)$ and $b_i$ neighbors of $y$ in $S_{i+1}(x).$ Evidently
$\Gamma$ is {\it regular} of valency $k=b_0,$ or {\it $k$-regular.}
A $k$-regular simple graph $\Gamma$ is {\it strongly regular} if
there exist integers $\lambda$ and $\mu$ such that  every adjacent
pair of vertices has $\lambda$ common neighbors, and every
nonadjacent pair of vertices has $\mu$ common neighbors.

The Hamming space $F_q^n$ consists of the $q^n$ vectors of length
$n$ over the alphabet $\{0,1,...,q-1\}, \ q\geq 2$. It is endowed
with the Hamming distance $d$ where $d(x,y)$ is the number of
coordinate positions in which $x$ and $y$ differ. It can be viewed
as a graph $L_n(q)$ with vertex set given by the vector space
$F_q^n$ (where $F_{q}$ is the field of $q$ elements) where $\{x,y\}$
is an edge of $L_n(q)$ iff $d(x,y)=1$. This Hamming graph is the
Cayley graph on the additive group $F_q^n$ when we take the
generator set $S=\{xe_{i}: x\in (F_{q})^{\times}, \, 1\leq i \leq n
\}$ where the $e_{i}=(0,...,0,1,0,...0)$ are the standard basis
vectors of $F_q^n$. It was shown in \cite{lev97,eros} that for any
$n\geq 2,$ $q\geq 2$ and $r\geq 1$,
\begin{equation}
\label{e2} N(L_n(q),r)=q\sum_{i=0}^{r-1}\left({n-1} \atop{i}\right)
(q-1)^i.
\end{equation}
For the particular case $n=2$ the Hamming graph $L_2(q)$ is the
{\it lattice graph} over $F_{q}$. This graph is strongly regular
with parameters $v=q^2,$ $k=2(q-1),$ $\lambda=q-2,$ $\mu=2,$ and
from (\ref{e2}) we get $N(L_2(q),1)=q$ and $N(L_2(q),2)=q^2.$

The Johnson graph $J_e^n$ is defined on the subset $V=J_e^n\subseteq
F_{2}^{n}$ consisting of all $\{0,1\}$-vectors with exactly $e$
entries equal to $1,$ for a fixed $1\leq e \leq n-1.$ On $J_e^n$ the
Johnson distance is defined as half the (even) Hamming distance, and two
vertices $x$, $y$ are joined by an edge iff they are at Johnson
distance $1$ from each other. In general $J_e^n$ is not a  Cayley
graph although the notion of errors being represented by edges makes
sense all the same. In particular, two vertices are at distance $1$
from each other iff one is obtained from the other by the
interchange of two coordinate positions. In \cite{lev97,eros} it was
shown that
\begin{equation}
\label{e3}N(J_e^n,r)=n\sum_{i=0}^{r-1}\left({e-1} \atop{i}\right)
\left ({n-e-1} \atop{i}\right) \frac 1{i+1}.
\end{equation}
for any $n\geq 2,\,e\geq 1$ and $r\geq 1.$ In the particular case
$e=2$ and $n \geq 4$ the Johnson graph $J_2^n$ is the {\it
triangular graph} $T(n)$. As vertices it has the 2-element subsets
of an $n$-set and two vertices are adjacent iff they are not
disjoint. This graph  is strongly regular with parameters
$v=\frac{n(n-1)}{2},$  $k=2(n-2),$  $\lambda=n-2,$ $\mu=4,$ and from
(\ref{e3}) we obtain $N(T(n),1)=n$ and $N(T(n),2)=\frac{n(n-1)}{2}.$

These two results were the first analytic formulas for the
reconstruction problem we are interested in. Their uniformity
depends on the fact that these graphs are distance--regular. What
then are the general results for simple graphs, regular graphs and
Cayley graphs? We start with a few observations from \cite{kls06}
for any connected simple graphs $\Gamma=(V,E)$. In the spirit of
distance regularity we put $k_i(x)=|S_i(x)|$ and define numbers
$c_i(x,y), \, b_i(x,y)$ and $a_i(x,y)$ for any two vertices $x\in V$
and $y\in S_i(x)$ such that
$$c_i(x,y)=|\{z\in S_{i-1}(x)\,\,:\,\, d(z,y)=1\}|\,,$$
$$b_i(x,y)=|\{z\in S_{i+1}(x)\,\,:\,\, d(z,y)=1\}|\,,$$
$$a_i(x,y)=|\{z\in S_i(x)\,\,: \,\,d(z,y)=1\}|\,.$$

From this $a_1(x,y)=a_{1}(y,x)$ is the number of triangles over the edge
$\{x,y\}$ and $c_2(x,y)$ is the number of common neighbors of $x\in V$
and $y\in S_2(x).$ Let
\begin{equation}
\label{e4} \lambda=\lambda(\Gamma)=\max_{x\in V,\ y\in
S_1(x)}a_1(x,y)
\end{equation}
\begin{equation}
\label{e5}\mu =\mu(\Gamma)=\max_{x\in V,\ y\in S_2(x)}c_2(x,y).
\end{equation}
Since $|B_r(x)\cap B_r(y)|>0$ for $x\neq y$ if and only if $1\le
d(x,y)\le 2r=d(\Gamma),$ we have
\begin{equation}
\label{e6} N(\Gamma,r)=\max_{1\le s \le 2r}N_s(\Gamma,r)
\end{equation}
where $N_s(\Gamma,r)=\max\{| B_r(x)\cap
B_r(y)|\,\,:\,\,x,y\in V;\ d(x,y)=s\}.$ In  particular, $N_1(\Gamma,1)=\lambda +2\,$ and
$\,N_2(\Gamma,1)=\mu$ so that
\begin{equation}
\label{e7} N(\Gamma,1)=\max(\lambda +2, \mu).
\end{equation}
One can easily check that using this formula for the lattice graph
$L_2(q)$ and the triangular graph $T(n)$ we obtain again the earlier
formulas (\ref{e2}) and (\ref{e3}). Indeed, since $\lambda=n-2$ and
$\mu=4$ for $T(n), n \geq 4,$ we have $N(T(n),1)=n$ from (\ref{e7}).
By the same reason we have $N(L_2(q),1)=q$ since $\lambda=q-2$ and
$\mu=2$ in this case.

We have no general results for $N(\Gamma,r)$ when $\Gamma$ is a
regular graph. The numbers $c_i(x,y)$ and $b_i(x,y)$ usually depend
on $y\in S_i(x)$ and this causes difficulties when searching for
general estimates of $N(\Gamma,r)$. However, some bounds on
$N(\Gamma,1)$ and $N(\Gamma,2)$ were obtained in \cite{kls06}. Here
it is assumed that $\Gamma$ is  connected, $k$-regular of diameter
$d(\Gamma)\ge 2$ with $v\ge 4$ vertices and parameters $0\le \lambda
\le k-2,$ $1\le \mu \le k,$ where $2\le k \le v-2.$ \vspace{2mm}
\begin{theorem} \label{th1} For any
$k$-regular graph $\Gamma$,
\begin{equation}
\label{e8} N(\Gamma,1)\le \frac 12 (v+\lambda).
\end{equation}
\end{theorem}
\vspace{2mm} This theorem is proved by checking that $\lambda+2\le
\frac 12 (v+\lambda)$ and $\mu \le \frac 12 (v+\lambda).$ The first
inequality takes place since $k\le v-2$ and $\lambda \le k-2.$
Moreover, there is equality only if $\lambda=v-4$ and $k=v-2$. The
second inequality is true since counting edges between $S_1(x)$ and
$S_2(x)$ for any $x\in V$ we have $\sum_{y\in
S_1(x)}(k-1-a_1(x,y))=\sum_{z\in S_2(x)}c_2(x,z).$ From (\ref{e4}), (\ref{e5}) and the fact that $k_2(x)\le v-k-1$ we get $
k(k-1-\lambda)\le \mu k_2(x)\le \mu (v-k-1)$ with equality if and
only if $\Gamma$ is strongly regular. Let us note here that the
equality $k(k-1-\lambda)=\mu (v-k-1)$ is well-known for strongly
regular graphs. From this and the fact that $1\le \mu\le k$ we have
$k-1-\lambda\le v-k-1$ and hence $\mu\le k\le \frac 12(v+\lambda)$
is valid for any regular graph $\Gamma$. By taking into account
these two inequalities for $\lambda$ and $\mu$ we get (\ref{e8})
from (\ref{e7}). Moreover, (\ref{e8}) is attained on the strongly
regular $t$-partite graph $K^{(t)}_{k-\lambda}$ with $t(k-\lambda)$
vertices partitioned into $t\geq 2$ parts, where $t=\frac
{2k-\lambda}{k-\lambda}$ is an integer, and with edges connecting
any two vertices of different parts.

\vspace{2mm}\begin{theorem} \label{th2} For any $k$-regular graph
$\Gamma$ we have
\begin{equation}
\label{e9} N_2(\Gamma,2)\ge\mu\left(k-1-\frac 34(\mu
-1)(N(\Gamma,1)-2)\right)+2.
\end{equation}
\end{theorem}

\vspace{2mm} In proving (\ref{e9}) the linear programming problem
arises for the vertex subset
$U=\bigcup_{i=1}^{\mu}B_1(z_i)\setminus\{x,y\},$ where $x,y\in V$
with $d(x,y)=2$ and $z_i, \ i=1,...,\mu,$ are the vertices at
distance 1 from both $x$ and $y$. The task is to minimize
$|U|=\sum_{h=1}^{\mu}u_h$ for nonnegative numbers $u_h$ satisfying
the following conditions
$$ \sum_{h=1}^{\mu}u_hh^2\ge \mu (k-1),$$
$$ \sum_{h=1}^{\mu}u_hh{h\choose 2}\le {\mu \choose 2}(N(G,1)-2),$$
where $u_h=|U(h)|/h,$ and $U(h)$ is the set of vertices in $U$
belonging to $h$ sets $B_1(z_i)$, $i=1,...,\mu$.

The details of the proofs for Theorems~1 and 2 can be found in \cite{kls06}. From the last  theorem one can immediately get the following corollaries.
\vspace{2mm}
\begin{corollary}
\label{cor1} For a $k$-regular graph $\Gamma,$
\newline (i) \,\, if $\mu=1,$ then $N_2(\Gamma,2)\ge k+1$;
\newline (ii) \ if $\mu=2$ and $N(\Gamma,1)=2$, then $N_2(\Gamma,2)\ge 2k$;
\newline (iii) if $\mu=3$ and $N(\Gamma,1)=3,$ then $N_2(\Gamma,2)\ge 3k-5.$
\end{corollary}
\vspace{2mm}
\begin{corollary}
\label{cor2} Let $\Gamma$ be a $k$-regular graph  without triangles
or pentagons, with $\mu \ge 2$ and $k\ge 1+\frac34(\mu -1)\mu.$ Then
\begin{equation}
\label{e10}N_2(\Gamma,2)\ge N_1(\Gamma,2). \end{equation}
\end{corollary}
\vspace{2mm}

Actually, since $\Gamma$ does not contain triangles or pentagons we have
$N_1(\Gamma,2)=2k$ and $N(\Gamma,1)=\mu$ by~(\ref{e7}) since
$\lambda =0$ and $\mu \ge 2.$ Using~(\ref{e9}) we get
$$N_2(\Gamma,2)-2k\ge (\mu-2)(k-1-\frac 34(\mu -1)\mu)\ge 0,$$
and finally we obtain (\ref{e10}).

In the remainder of this section it is assumed that
$\Gamma=Cay(G,S)$ is a Cayley graph on the group $G$ for the
generator set $S$. Let us put $S^{0}=\{e\}$ and set
$S^{i}=SS^{i-1}$. Moreover, by vertex--transitivity it is sufficient
to consider only the spheres and balls with center $e$ so that
$S_{i}=S_{i}(e)$.

\vspace{2mm}
\begin{lemma} \label{lem1}  For any Cayley
graph $\Gamma$ on the group $G$ and for $i>0$ we have
$S_{i}=S^{i}\setminus (S^{i-1}\cup S^{i-2}\cup...\cup S^{0})$. In
particular, $\mu$ is the maximum number of representations of an
element in $S^{2}\setminus (S\cup S^{0})$ as a product of two elements
of $S$ and $\lambda$ is the maximum number of representations of an
element in $S$ as a product of two elements of $S$, i,e.
$$\lambda(\Gamma) =\max_{s\in S}\mid \{(s_is_j)\in S^2 \,:\,s=s_is_j
\} \mid,$$ $$\mu(\Gamma) =\max_{s\in S^2\setminus (S\cup S^{0})}
\mid \{(s_is_j)\in S^2 \,:\,s=s_is_j \} \mid.$$
\end{lemma} \vspace{2mm}

\noindent This lemma allows us to find $N(\Gamma,1)$ from (\ref{e7})
for a general Cayley graph. The results for
estimating the values $N(\Gamma,r)$ for small $r$ in  Cayley
graphs on the symmetric group $Sym_n$ will be presented in the next section
when the generator set $S$ consists of transpositions.

\section{The reconstruction of permutations in Cayley graphs generated by transpositions}

\noindent Let $Sym_n$ be the symmetric group on $n$ symbols. We write a
permutation $\pi$ in one--line notation as
$\pi=[\pi_1,\pi_2,\ldots,\pi_n]$ where $\pi_i=\pi(i)$ for every
$i\in\{1,\ldots,n\}.$

For the {\it transposition Cayley graph} $Sym_n(T)$ on $Sym_n$ the
generator set consists of all transpositions  $T=\{t_{i,j}\in Sym_n,
\ 1 \leq i < j \leq n\}, \ |T|=\left({n}\atop{2}\right),$ where
$t_{i,j}$ interchanges positions $i$ and $j$ when multiplied on  the
right, i.e., $[\ldots,\pi_i,\ldots,\pi_j,\ldots]\cdot
t_{i,j}=[\ldots, \pi_j,\ldots,\pi_i,\ldots].$ For $x,\,y\in Sym_{n}$
the distance $d(x,y)$ is the least number of transpositions
$t_{1},\,...,t_{r}$ such that $x\cdot t_{1}\cdot...\cdot t_{r}=y,$
or $ t_{1}\cdot...\cdot t_{r}=x^{-1}\cdot y$. As any $k$-cycle can
be written as a product of $k-1$ transpositions (but no fewer), the
diameter of $Sym_n(T)$ is $(n-1).$ The graph is bipartite since any
edge joins an even permutation to an odd permutation. The symmetry
properties of $Sym_n(T)$ have been discussed in \cite{LJD93}. The
graph  is edge--transitive but not distance--regular and hence not
distance--transitive. All these properties and other basic facts are
collected in the following statements.

\vspace{2mm}
\begin{lemma}
\label{lem2} The transposition graph $Sym_n(T), \  n \geq 3,$
\newline (i) is a connected bipartite
$\left({n}\atop{2}\right)$-regular graph of order $n!$ and diameter
$(n-1)$; \newline (ii) is not  distance--regular and hence not
distance--transitive; \newline (iii) it does not contain
subgraphs isomorphic to $K_{2,4},$ and each of its vertices belongs to
$\left({n} \atop{3}\right)$ subgraphs isomorphic to $K_{3,3}.$
\end{lemma}

\vspace{2mm}\noindent (Here $K_{p,q}$ is the complete bipartite graph
with $p$ and $q$ vertices in the two parts, respectively.)

\vspace{2mm}
\begin{theorem} \label{th3} For any \ $n\ge 3$ we have
$N(Sym_n(T),1)=3.$
\end{theorem}
\vspace{2mm} This means that any unknown permutation is uniquely
reconstructible from 4 distinct permutations at transposition
distance at most one from the unknown permutation. The proof of
these statements is based on considering a permutation $\pi \in
Sym_n$ in cycle notation, with {\it cycle type} ${\rm
ct}(\pi)=1^{h_{1}}2^{h_{2}}..\,n^{h_{n}}$, where $h_{i}$ is the
number of cycles of length $i.$ In particular
$\sum_{i}^{n}ih_{i}=n$. The permutation $\pi$ can be also presented
as a product of a least number of transpositions. Each such product
represents a shortest path in $Sym_n(T)$ from $e$ to $\pi$. The
number of such paths was obtained in \cite{den}. This result is
based on Ore's theorem on the number of trees with $n$ labeled
vertices and presented by the following theorem. \vspace{2mm}
\begin{theorem}\cite{den}
\label{th4} Let $\pi \in Sym_n$ have cycle type ${\rm
ct}(\pi)=1^{h_{1}}2^{h_{2}}...n^{h_{n}}$, consisting of
$\sum_{j=1}^{n}\,h_{j}=n-i$ cycles where $1\leq i \leq n-1.$ Then
the number of distinct ways to express $\pi$ as a product of $i$
transpositions is equal to
$$i!\prod_{j=1}^n\left(\frac{j^{j-2}}{(j-1)!}\right)^{h_j}.$$
\end{theorem}
\vspace{2mm} According to the above theorem, the following lemma
gives us formulas for the numbers $c_i(\pi):=c_i(\pi,e),$
$b_i(\pi):=b_i(\pi,e),$ $a_i(\pi)=a_i(\pi,e), 1\leq i \leq n-1$
where $e$ is the identity element of the transposition Cayley graph
$Sym_n(T).$ \vspace{2mm}
\begin{lemma} \label{lem3} In the transposition
graph $Sym_n(T)$ the sets $S_{i}=S_{i}(e), 1\leq i \leq n-1,$ are
the permutations consisting of $(n-i)$ disjoint cycles, counting
also $1$-cycles. For any $\pi \in S_i$ with cycle type ${\rm
ct}(\pi)=1^{h_{1}}2^{h_{2}}...\,n^{h_{n}},$ we have $a_i(\pi)=0$ and
$$ c_i(\pi)=\frac 12 \left(\sum_{j=1}^nj^2h_j-n\right), \
b_i(\pi)=\frac 12 \left(n^2-\sum_{j=1}^nj^2h_j\right).$$
\end{lemma}
\vspace{2mm} In particular, since $a_i(\pi)=0$ for any $1\leq i \leq
n-1,$ then from this lemma and by~(\ref{e4}) we have
$\lambda(Sym_n(T))=0.$ Moreover, it is well-known that two
permutations are conjugate by an element of $Sym_n$ if and only if
they have the same cycle type. If
$(1^{h_{1}}2^{h_{2}}...\,n^{h_{n}})^{G}$ denotes the conjugacy class
of an element of cycle type $1^{h_{1}}2^{h_{2}}...\,n^{h_{n}}$ then
it is shown in \cite{kls06} that  $S_{i}, \, 1\leq i \leq n-1,$ is
the disjoint union
\begin{equation}
\label{e11} S_i=\bigcup_{h_1+h_2+\cdots +h_n=n-i}
(1^{h_{1}}\,2^{h_{2}}\,...\,n^{h_{n}})^{G},
\end{equation}where
\begin{equation}
\label{e12} |(1^{h_{1}}2^{h_{2}}...\,n^{h_{n}})^{G}|= \frac
{n!}{1^{h_1}h_1!2^{h_2}h_2!\cdots n^{h_n}h_n!}.
\end{equation}

Hence, from~(\ref{e11}) we have $S_{2}= (1^{n-3}\,3^{1})^{G}\,\cup
(1^{n-4}\,2^{2})^{G}$ and then by Lemma~\ref{lem3} we get
$c_2(\pi)=3$ if $ct(\pi)=1^{n-3}\,3^1,$ and $c_2(\pi)=2$ if
$ct(\pi)=1^{n-4}\,2^2.$ From these and~(\ref{e5}) we have
$\mu(Sym_n(T))=3,$ and therefore, by~(\ref{e7}) we get
Theorem~\ref{th3}. Moreover, there are no subgraphs isomorphic to
$K_{2,4}$ in $Sym_n(T)$ since $\mu(Sym_n(T))=3.$ The number
$\left({n} \atop{3}\right)$ of subgraphs isomorphic to $K_{3,3}$ and
having $e$ as one of its vertices is obtained from~(\ref{e12}) for
any $\pi\in (1^{n-3}\,3^1)^G.$ By vertex--transitivity the same
holds for any vertex in $Sym_n(T)$ (see condition (iii) in
Lemma~\ref{lem2}).

So, any unknown permutation is uniquely reconstructible from $4$
distinct permutations at transposition distance at most $1$ from the
unknown permutation. As the following shows, in the case of at most
two transposition errors
the reconstruction of the permutation $\pi$ requires many more  distinct
$2$--neighbors of $\pi$.

\vspace{2mm}
\begin{theorem}
\label{th5} For $n\ge 3$ we have
\begin{equation}
\label{e13} N(Sym_n(T),2)=\frac 32(n-2)(n+1).
\end{equation}
\end{theorem}

\vspace{2mm} \noindent This  follows from the fact that the
normalizer of $T$ is $G=Sym_n$ itself and from the following lemma.
\vspace{2mm}
\begin{lemma} \label{lem4} For any $\pi\in S_i, \,
1\leq i \leq n-1$ the number of  vertices in
$(1^{h_{1}}2^{h_{2}}...\,n^{h_{n}})^{G}$ at a given distance from
$\pi$ depends only on the conjugacy class to which $\pi$ belongs.
\end{lemma}
\vspace{2mm}

To prove Theorem~\ref{th5} it is therefore sufficient to consider
the numbers of vertices in all subsets of $B_2(e)$ at minimal
distance at most 2 from a given vertex $\pi \in S_i, \, 1\leq i \leq
4.$ By~(\ref{e11}) we have $S_1= (1^{n-2}\,2^{1})^{G}$, $S_2=
(1^{n-3}\,3^{1})^{G}\,\cup (1^{n-4}\,2^{2})^{G}$,
$S_3=(1^{n-4}\,4^1)^G\cup\,(1^{n-5}\,2^1\,3^1)^G\cup\,(1^{n-6}\,2^3)^G,$
$S_4=(1^{n-5}\,5^1)^G\cup\,(1^{n-6}\,2^1\,4^1)^G\cup\,(1^{n-6}\,3^2)^G\bigcup\,(1^{n-7}\,2^2\,3^1)^G\bigcup\,$
$(1^{n-8}\,2^4)^G.$ It is shown that $N_4(Sym_n(T),2)=20$ for $n\ge
5,$ $N_3(Sym_n(T),2)=12$ for $n\ge 4,$ $N_2(Sym_n(T),2)=\frac 32
(n-2)(n+1)$ and $N_1(Sym_n(T),2)=n(n-1)$ for all $n\ge 3.$ From
these and by~(\ref{e6}) one can conclude~(\ref{e13}).

The statements of Theorem~\ref{th5} and Corollary~\ref{cor2} are
generalized in the following conjecture. \vspace{2mm}
\begin{conjecture}
\label{con1} For any $\pi\in (1^{n-3}\,3^1)^G,$ for any $r\ge 1$ and
$n\ge 2r+1$ we have $$N(Sym_n(T),r)=N_2(Sym_n(T),2)=|B_r(I)\bigcap
B_r(\pi)|.$$
\end{conjecture}
\vspace{2mm} Now let us consider the {\it bubble--sort graph}
$Sym_n(t).$ This is the Cayley graph on $Sym_n$ for the generator
set $t=\{t_{i,i+1}\in Sym_n, \ 1 \leq i < n\}, \ |t|=n-1.$ These
{\it bubble--sort transpositions} are  $2$-cycles \ $t_{i,i+1}$
interchanging  $i$ and $i+1$ and determine the graph distance in
$Sym_n(t)$ in the usual way. It is known that the
diameter of $Sym_n(t)$ is $\left({n}\atop{2}\right).$ 
\vspace{2mm}
\begin{lemma}
\label{lem5}  The bubble--sort graph $Sym_n(t), \ n \geq 3,$
\newline (i) is a connected bipartite $(n-1)$-regular graph of order $n!$
and diameter $\left({n}\atop{2}\right);$ \newline (ii)  it does not
contain subgraphs isomorphic to $K_{2,3};$
\newline (iii) each of its vertices belongs to $\left({n-2} \atop{2}\right), n \geq 4,$ subgraphs
isomorphic to $K_{2,2}.$
\end{lemma}
\vspace{2mm} The symmetry properties of the bubble--sort graph were
discussed in \cite{LJD93} where it was shown that this graph is not
distance--regular. As it is bipartite there are  no triangles and
hence $\lambda(Sym_n(t))=0.$
If an element $\pi\in S_{2}(e)$ has at least two neighbors
$t_{i,i+1}\neq t_{j,j+1}$ in  $S_{1}(e)$ then  necessarily
$t_{i,i+1}t_{j,j+1}=\pi=t_{j,j+1}t_{i,i+1}$ with $\{j,j+1\}$ and
$\{i,i+1\}$ disjoint. It suffices to verify this for permutations on
$4$ letters. Hence there are at most two such neighbors and so
$\mu(Sym_n(t))=2. $ It can be also verified that we have
$N_4(Sym_n(t),2)=4$ for $n\ge 5,$ $N_3(Sym_n(t),2)=2$ for $n\ge 4,$
$N_2(Sym_n(t),2)=N_1(Sym_n(t),2)=2(n-1)$ for $n\ge 3.$ From all
these and by~(\ref{e6}) and~(\ref{e7}) we get the following theorem.
\vspace{2mm}
\begin{theorem}  \label{th6} For any $n\ge 3$ we have
$$N(Sym_n(t),1)=2\quad{\rm and} \quad N(Sym_n(t),2)=2(n-1).$$ 
\end{theorem}
\vspace{4mm} Almost the same results appear for the {\it star Cayley
graph} $Sym_n(st)$ generated by the set of {\it
prefix--transpositions} $st=\{(1,i)\in Sym_n, \ 1 < i \leq n\}, \
|st|=n-1.$ It is one of the most investigated graphs in the theory
of interconnection networks since many parallel algorithms are
efficiently mapped on the star Cayley graph.

\vspace{2mm} \begin{lemma} \cite{AK89} \label{lem6} The star Cayley
graph $Sym_n(st), n \geq 3,$ is a connected bipartite
$(n-1)$-regular graph of order $n!$ with diameter $\lfloor
\frac{3(n-1)}{2}\rfloor$.
\end{lemma}
\vspace{2mm}

The star Cayley graph $Sym_n(st)$ is not distance--regular for $n\ge
4$ \cite{LJD93} and has no cycles of lengths of 3, 4, 5 or 7. Hence
$\lambda(Sym_n(st))=0$ and $\mu(Sym_n(st))=1$. Moreover, it is easy
to verify that  $N_4(Sym_n(st),2)=4$ for $n\ge 5,$
$N_3(Sym_n(st),2)=4$ for $n\ge 4,$ $N_2(Sym_n(st),2)=2(n-1)$ for
$n\ge 5$ and $N_1(Sym_n(st),2)=2(n-1)$ for $n\ge 4.$ From these
properties and by~(\ref{e6}) and~(\ref{e7}) we get the following
theorem.

\vspace{2mm}
\begin{theorem} \label{th7} For any $n\ge 4$ we have
$$N(Sym_n(st),1)=2 \quad{\rm and} \quad N(Sym_n(st),2)=2(n-1).$$ 
\end{theorem}
\vspace{2mm}

Thus, in the bubble--sort and star Cayley graphs any unknown
permutation $\pi$ is uniquely reconstructible from 3  distinct
$1$--neighbors of $\pi$. 
Similarly, for the unique reconstruction of $\pi$ from neighbors at
distance at most $2$ we see that any $2n-1$ distinct $2$--neighbors
of $\pi$ are sufficient.  These two graphs are examples for which
the inequality $(ii)$ in Corollary~(\ref{cor1}) is attained.

\vspace{5mm}
\section*{Acknowledgment}
\noindent The research was partially supported by the RFBR grant
06--01--00694.

\end{document}